\documentclass{amsart}

\usepackage{amssymb,amsmath,amsfonts,amsthm, mathrsfs}
\usepackage[all,dvips]{xy}
\newtheorem{theorem}{Theorem}[section]
\newtheorem{thm}[theorem]{Theorem}
\newtheorem{lem}[theorem]{Lemma}
\newtheorem*{cor}{Corollary}

\theoremstyle{remark}
\newtheorem*{remark}{Remark}
\newtheorem*{definition}{Definition}
\theoremstyle{remark}

\newcommand{\onto}{\twoheadrightarrow}

\newcommand{\tor}{\mathop{\mathrm{Tor}}\nolimits}

\newcommand{\x}{\mathbf x}

\renewcommand{\l}{\ell}

\renewcommand{\i}[1]{\mathfrak{#1}}

\begin{document}
\title[CHARACTERIZATIONS OF REGULAR LOCAL RINGS] {CHARACTERIZATIONS OF REGULAR LOCAL RINGS\\ IN
POSITIVE CHARACTERISTICS}
\author{Jinjia Li}
\address{Department of Mathematics, Syracuse University, 215 Carnegie, Syracuse, NY 13244}
\curraddr{} \email{jli32@syr.edu}


\date{\today}


\begin{abstract}
In this note, we provide several characterizations of regular
local rings in positive characteristics, in terms of the
Hilbert-Kunz multiplicity and its higher $\tor$ counterparts $\i
t_i=\underset{n \to \infty}{\lim} \l(\tor_i(k,{}^{f^n}\!\!
R))/p^{nd}$. We also apply the characterizations to improve a
recent result by Bridgeland and Iyengar in the characteristic $p$
case. Our proof avoids using the existence of big Cohen-Macaulay
modules, which is the major tool in the proof of Bridgeland and
Iyengar.
\end{abstract}

\maketitle

\section{Introduction}

Let $(R,\i m,k)$ be a $d$-dimensional local ring of characteristic
$p>0$. The Frobenius endomorphism $f_R:R \to R$ is defined by
$f_R(r)=r^p$ for $r \in R$. Each iteration $f_R^n$ defines a new
$R$-module structure on $R$, denoted ${}^{f^n}\!\! R$, for which
$a\cdot b=a^{p^n}b$. For any $R$-module $M$, $F_R^n(M)$ stands for
$M\otimes_R {}^{f^n}\!\! R$, the $R$-module structure of which is
given by the base change along the Frobenius endomorphism. When
$M$ is a cyclic module $R/I$, it is easy to show that $F^n(R/I)
\cong R/I^{[p^n]}$, where $I^{[p^n]}$ denotes the ideal generated
by the $p^n$-th power of the generators of $I$.

In what follows, $\l(-)$ denotes the length function.

For any $\i m$-primary ideal $I$, the Hilbert-Kunz multiplicity of
$R$ with respect to $I$ was first introduced by Monsky in
\cite{Mo}:
\[e_{\text{HK}}(I,R) = \underset{n\to \infty} {\lim}
\l(F^n(R/I))/p^{nd}.\] The Hilbert-Kunz multiplicity of $R$ is
$e_{\text{HK}}(R)=e_{\text{HK}}(\i m,R)$. We also frequently write
$e_{\text{HK}}(I)=e_{\text{HK}}(I,R)$. It has been shown by many
authors that the Hilbert-Kunz multiplicity encodes subtle
information about the singularity of $R$. One such example is the
following characterization of the regularity due to Watanabe and
Yoshida; see \cite{W-Y}, \cite{H-Y}.

\begin{thm}
If $R$ is unmixed, then it is regular if and only if
$e_{\text{HK}}(R)=1$.
\end{thm}

It is natural to ask whether the higher $\tor$ counterparts of the
Hilbert-Kunz multiplicity, which are defined below, can encode
similar information on the singularity of the ring.
\begin{definition}Let $R$ be a $d$-dimensional local ring of characteristic $p>0$.
Let $I$ be any $\i m$-primary ideal. Define
\[\i t_i(I,R)= \underset{n \to \infty}{\lim} \l(\tor_i(R/I,{}^{f^n}\!\!
R))/p^{nd}.\]
\end{definition}
Seibert has shown that such limits always exist \cite{Se}. In the
sequel, we also write $\i t_i(R)=\i t_i(\i m,R)$. The main result
of this note is the following theorem in Section 2:

\medskip

\noindent \textbf{Main Theorem} (Theorem~\ref{main}). Let $(R, \i
m, k)$ be a $d$-dimensional local ring of characteristic $p>0$.
Then the following are equivalent:
\begin{enumerate}
  \item[(i)] $R$ is regular,
  \item[(ii)] $\i t_1(R)=0$,
  \item[(iii)] $\i t_2(R)=0$,
  \item[(iv)] $e_{\text{HK}}(R)- 1 = \i t_1(R).$
\end{enumerate}

The proof of this theorem is inspired by the work of Huneke and
Yao \cite [2.1]{H-Y} and Blickle and Enescu \cite{B-E}.

\medskip

In Section 3, we apply our main theorem to slightly generalize the
positive characteristic case of a recent result by Bridgeland and
Iyengar \cite [1.1] {B-I}. The result of Bridgeland and Iyengar
states that if $R$ contains a field and if $k$ is a direct summand
of $H_0(C_\bullet)$, where $C_\bullet$ is a perfect complex of
length exactly $d$ such that $\l(H_i(C_\bullet))<\infty$ for all
positive $i$, then $R$ must be regular. Using Theorem~\ref{main},
we are able to show that in the positive characteristic case, not
only $k$, but also the first syzygy of $k$ cannot be a direct
summand of $H_0(C_\bullet)$ for such $C_\bullet$ unless $R$ is
regular.

Our proof of this result is quite different from that of
Bridgeland and Iyengar. In particular, we avoid using the
existence of big Cohen-Macaulay modules.

\section{The Main Result}
The following fact plays a crucial role in this paper. It is
contained in \cite [1.1] {D}. We provide a sketch of the proof
here for the completeness of this paper.
\begin{lem}
\label{lemma} Let $(R, \i m, k)$ be a $d$-dimensional local ring
of characteristic $p>0$. Let $I=(\x)$ be an ideal generated by a
system of parameters $\x=x_1, ..., x_d$. Then $\i t_1(I, R)=0$.
\end{lem}
\begin{proof}
There is a surjection
\[H_1(\x^{[p^n]};R) \onto \tor_1(R/I,{}^{f^n}\!\!
R)\] and it is well known (see \cite [7.3.5] {R} or \cite [1.1]
{D}) that
\[\underset{n \to \infty}{\lim} \l(H_i(\x^{[p^n]};R))/p^{nd}=0, \text{ for }i>0.\]
\end{proof}

The following is the main theorem of this paper.
\begin{thm} \label{main} Let $(R, \i m, k)$ be a $d$-dimensional local ring of characteristic $p>0$. Then
\begin{enumerate}
  \item[(a)] $e_{\text{HK}}(R)- 1 \leq \i t_1(R)$;
  \item[(b)] $ \i t_1(R) - e_{\text{HK}}(R)+ 1 \leq \i t_2(R)$.
\end{enumerate}
Moreover, the following are equivalent:
\begin{enumerate}
  \item[(i)] $R$ is regular,
  \item[(ii)] $\i t_1(R)=0$,
  \item[(iii)] $\i t_2(R)=0$,
  \item[(iv)] $e_{\text{HK}}(R)- 1 = \i t_1(R).$
\end{enumerate}

\end{thm}

\begin{proof} We first prove (a) and (b).
Let $q$ be a power of $p$. Let $I$ be an $\i m$-primary ideal
generated by a system of parameters of $R$. Consider the following
filtration:
\[
0 \to Q_1 \to R/I^{[q]}\to k \to 0,
\]
\[
0\to Q_2\to Q_1\to k\to 0,
\]
\[
\vdots
\]
\[
0\to k\to Q_l\to k\to 0.
\]
Applying $-\otimes{}^{f^n}\!\! R$ to the above short exact
sequences, we obtain the long exact sequences
\[
\xymatrix@C=15pt@R=26pt{
\cdots \ar[r] & \tor_2(k,{}^{f^n}\!\! R) \ar[r]& \tor_1(Q_1, {}^{f^n}\!\! R)\ar[r] & \tor_1(R/I^{[q]}, {}^{f^n}\!\! R) \ar `[r] `[l] `[dlll] `[dll] [dll] & \\
 &\tor_1(k,{}^{f^n}\!\! R) \ar[r] & F^n(Q_1) \ar[r]& F^n(R/I^{[q]}) \ar[r] & F^n(k) \ar[r] & 0,\\
}
\]
\[
\xymatrix@C=15pt@R=26pt{
\cdots \ar[r] & \tor_2(k,{}^{f^n}\!\! R) \ar[r]& \tor_1(Q_2, {}^{f^n}\!\! R)\ar[r] & \tor_1(Q_1, {}^{f^n}\!\! R) \ar `[r] `[l] `[dlll] `[dll] [dll] & \\
 &\tor_1(k,{}^{f^n}\!\! R) \ar[r] & F^n(Q_2) \ar[r]& F^n(Q_1) \ar[r] & F^n(k) \ar[r] & 0,\\
}
\]
\[ \vdots \]
\[
\xymatrix@C=15pt@R=26pt{
\cdots \ar[r] & \tor_2(k,{}^{f^n}\!\! R) \ar[r]& \tor_1(k, {}^{f^n}\!\! R)\ar[r] & \tor_1(Q_l, {}^{f^n}\!\! R) \ar `[r] `[l] `[dlll] `[dll] [dll] & \\
 &\tor_1(k,{}^{f^n}\!\! R) \ar[r] & F^n(k) \ar[r]& F^n(Q_l) \ar[r] & F^n(k) \ar[r] & 0.\\
}
\]

It follows that \begin{equation} (\l (R/I^{[q]})-1)\cdot
\l(\tor_1(k,{}^{f^n}\!\! R))+  \l(F^n(R/I^{[q]})) \geq
\l(R/I^{[q]})\cdot \l(F^n(k))
\end{equation}
and
\begin{align*}
(\l (R/I^{[q]})-1)\cdot \l(\tor_2(k,{}^{f^n}\!\! R))+&
\l(\tor_1(R/I^{[q]},{}^{f^n}\!\! R))+ \l(R/I^{[q]})\cdot
\l(F^n(k))\\
&\geq \l(R/I^{[q]})\cdot \l(\tor_1(k, {}^{f^n}\!\! R))+
\l(F^n(R/I^{[q]})) \tag{2}.
\end{align*}

Divide both sides of (1) by $p^{nd}$ and let $n \to \infty$ to
obtain
\[(\l (R/I^{[q]})-1)\cdot \i t_1(R) + q^d \cdot e_{\text{HK}}(I) \geq \l(R/I^{[q]})e_{\text{HK}}(R) \tag{3}.\]

Dividing (3) by $q^d$ and letting $q \to \infty$ then yields
\[e_\text{HK}(I)\cdot \i t_1(R)+e_\text{HK}(I) \geq e_{\text{HK}}(I)\cdot e_{\text{HK}}(R). \]
Hence \[e_{\text{HK}}(R) \leq \ 1 + \i t_1(R).\]

Similarly, from inequality (2) we can get (we need to apply
Lemma~\ref{lemma} here)
\[(\l (R/I^{[q]})-1)\cdot \i t_2(R) + \l(R/I^{[q]})e_{\text{HK}}(R) \geq \l(R/I^{[q]})\cdot \i t_1(R) + q^d \cdot e_{\text{HK}}(I)\tag{4}.\]
Therefore
\[ \i t_1(R) - e_{\text{HK}}(R)+ 1 \leq \i t_2(R).\]

For the second part of the theorem, it is clear that (i)
$\Rightarrow$ (ii), (iii) and (iv) due to the exactness of
Frobenius. We now prove (ii) $\Rightarrow$ (i) and (iii)
$\Rightarrow$ (iv) $\Rightarrow$ (ii).

\medskip

(ii)$\Rightarrow$(i). Note that inequality (3) is valid for any
$\i m$-primary ideal $I$ (not just for ideals generated by a
system of parameters). Since $\i t_1(R)=0$, inequality (3) becomes
the equality
\[q^d \cdot e_{\text{HK}}(I) = \l(R/I^{[q]})e_{\text{HK}}(R). \]
Taking $I=\i m$ immediately gives $q^d = \l(R/\i m^{[q]})$, which
forces $R$ to be regular by Kunz's Theorem \cite{K}.

\medskip

(iii)$\Rightarrow$(iv). Since $\i t_2(R)=0$, inequality (4)
becomes the equality
\[\l(R/I^{[q]})e_{\text{HK}}(R) = \l(R/I^{[q]})\cdot
\i t_1(R) + q^d \cdot e_{\text{HK}}(I).\] We therefore obtain (iv)
by dividing both sides by $q^d$ and taking the limits.

\medskip

(iv)$\Rightarrow$(ii). We can make a flat extension of $R$ to
assume that $k$ is infinite without changing any of the relevant
lengths. Let $I$ be a minimal reduction of $\i m$ which is
generated by a system of parameters of $R$. In this case, it is
well known (see, for instance, \cite[14.12]{Mat}) that the
Hilbert-Kunz multiplicity $e_\text{HK}(I)$ coincides with the
Hilbert-Samuel multiplicity $e(R)$. Taking $q=1$ in (3), we have
\[(\l (R/I)-1)\cdot \i t_1(R) +  e_{\text{HK}}(I) \geq \l(R/I)e_{\text{HK}}(R) \tag{5}.\]
Replacing $e_{\text{HK}}(R)$ by $1 + \i t_1(R)$ in (5), we get
\[e(R) \geq \l(R/I)+\i t_1(R) \geq \l(R/I). \tag{6}\]
On the other hand, since $I$ is a minimal reduction of $\i m$,
$e(R) \leq \l(R/I)$. This forces all the inequalities in (6) to be
equalities. Therefore $\i t_1(R)=0$.
\end{proof}

The inequalities (a) and (b) in Theorem~\ref{main} are far from
being the best possible bounds for the Hilbert-Kunz multiplicity.
For example, the following corollary gives better bounds when $R$
is Cohen-Macaulay.

\begin{cor} Let $(R, \i m)$ be a Cohen-Macaulay local ring of
characteristic $p>0$ and let $e=e(R)$ be the Hilbert-Samuel
multiplicity of $R$. Then
\[e_{\text{HK}}(R)-1 \leq \bigg(\dfrac{e-1}{e}\bigg)\i t_1(R).\]
\end{cor}

This follows easily from inequality (3) in the proof of
Theorem~\ref{main}. One can again assume the residue field of $R$
is infinite so that there is a minimal reduction $I$ of $\i m$
which is generated by a system of parameters. Since $R$ is
Cohen-Macaulay, we have
\[e=\l(R/I)=e_{\text{HK}}(I).\]
Take $q=1$ in (3) and then replace $\l(R/I)$ and
$e_{\text{HK}}(I)$ in (3) by $e$. We obtain
\[(e-1)\cdot \i t_1(R) +  e \geq e \cdot e_{\text{HK}}(R),\]
which gives the desired inequality.

\medskip

\begin{remark}
When $R$ is Cohen-Macaulay, we can argue exactly the same way as
in the proof of Theorem~\ref{main} to trivially generalize
Theorem~\ref{main} to cases of higher $\tor$. Namely, we have
\[\i t_i(R)-\i t_{i-1}(R)+\cdots+(-1)^{i-1}\i t_1(R) +
(-1)^ie_{\text{HK}}(R)+(-1)^{i+1}\geq 0 \text{  for all }i\geq 1\]
and $R$ being regular can be characterized by either the above
inequalities taking ``='' for some $i \geq 1$, or $\i t_i(R)$
being zero for some $i \geq 1$. However, the author does not know
if this generality can be true without the Cohen-Macaulay
assumption. The main obstruction here is, when $R$ is not
Cohen-Macaulay, we no longer have the higher $\tor$ (for $i\geq
2$) analog of Lemma~\ref{lemma}.
\end{remark}

\section{An improvement of Bridgeland-Iyengar's Result}
Recently, Bridgeland and Iyengar \cite [1.1]{B-I} proved the
following characterization for regular local rings.
\begin{thm} [Bridgeland-Iyengar] \label{bi} Let $(R,\i m,k)$ be a $d$-dimensional
local ring containing a field or of dimension $\leq3$. Assume
$C_\bullet$ is a complex of free $R$-modules with $C_i=0$ for $i
\notin [0, d]$, the $R$-module $H_0(C_\bullet)$ is finitely
generated, and $\l(H_i(C_\bullet))<\infty$ for $i>0$. If $k$ is a
direct summand of $H_0(C_\bullet)$, then $R$ is regular.
\end{thm}
Their proof of Theorem 3.1 uses the existence of balanced big
Cohen-Macaulay modules. Here we can apply Theorem~\ref{main} to
give a more direct proof in the positive characteristic case that
avoids using the existence of big Cohen-Macaulay modules.
Moreover, our proof also yields the same conclusion if (instead of
$k$) the first syzygy module of $k$ is a direct summand of
$H_0(C_\bullet)$.
\begin{thm} \label{app} Let $(R,\i m,k)$ be a $d$-dimensional local ring of characteristic $p>0$, $C_\bullet$ a complex of free
$R$-modules with $C_i=0$ for $i \notin [0, d]$, the $R$-module
$H_0(C_\bullet)$ finitely generated, and
$\l(H_i(C_\bullet))<\infty$ for $i>0$. If either $k$ or the first
syzygy of $k$ is a direct summand of $H_0(C_\bullet)$, then $R$ is
regular.
\end{thm}
\begin{proof}
By the same argument as in \cite[Lemma 2.2]{B-I}, we have a
surjection
\[H_1(F^n(C_\bullet)) \onto \tor_1(H_0(C_\bullet), {}^{f^n}\!\!
R)\]

It is well known that $\underset{n \to \infty}{\lim}
\l(H_1(F^n(C_\bullet)))/p^{nd}=0$; see \cite [1.7] {D}. So we are
done by Theorem~\ref{main}.

\end{proof}

\begin{remark}
Theorem~\ref{app} is still valid for rings containing a field or
of dimension $\leq 3$ (the exact same hypothesis on $R$ as in
Theorem ~\ref{bi}) although the proof requires the use of big
Cohen-Macaulay modules. To see this, one needs to use a result of
Schoutens \cite[Proposition 2.5]{Sc} to modify the original proof
of Bridgeland and Iyengar (their proof of \cite[Theorem 2.4]{B-I})
slightly. We leave the details here to the readers.
\end{remark}

It seems that the mixed characteristic case of the above result
remains unknown.

\specialsection*{ACKNOWLEDGEMENT} I wish to thank Graham Leuschke
and Claudia Miller for introducing the paper of Bridgeland and
Iyengar to me and useful discussions, and to thank C-Y. Jean Chan,
who had carefully read an earlier version of this paper and
suggested many improvements.

\bibliographystyle{amsalpha}

\begin{thebibliography}{10}
\bibitem[B-E] {B-E} M.Blickle and F. Enescu,
\emph{On rings with small Hilbert-Kunz multiplicity} Proc. Amer.
Math. Soc., 132, 2004,  no. 9, 2505--2509.
\bibitem[B-I] {B-I} T. Bridgeland and S. Iyengar, \textit{A criterion for regularity of local rings}, C. R. Math. Acad. Sci. Paris, 342, 2006, no. 10,
723--726.
\bibitem [D]{D} S. P. Dutta, \textit{Ext and Frobenius}, J.
Algebra, 127, 1989, 163--177.
\bibitem[H-Y] {H-Y} C. Huneke and Y. Yao, \textit{Unmixed local rings with minimal
Hilbert-Kunz multiplicity are regular}, Proc. Amer. Math. Soc.,
130, 2002, no. 3, 661--665.
\bibitem[K] {K} E. Kunz, \textit{Characterization of regular local
rings for characteristic p}, Amer. J. Math., 91 1969, 772-784.
\bibitem[Ma] {Mat} H.~Matsumura, \emph{Commutative Ring Theory}, Cambridge Stud. Adv. Math. 8,
Cambridge Univ. Press, Cambridge, 1986.
\bibitem[Mo] {Mo} P.~Monsky,
\emph{The Hilbert-Kunz function}, Math. Ann., 263, 1983, 43--49.
\bibitem[R] {R} P. Roberts, \textit{Multiplicities and Chern Classes in Local Algebra}, Cambridge University Press, 1998.
\bibitem [Sc]{Sc} H. Schoutens, \textit{On the vanishing of Tor of the absolute integral closure}, J. Algebra, 275, 2004,
567--574.
\bibitem [Se]{Se} G. Seibert, \textit{Complexes with homology of
finite length and Frobenius functors}, J. Algebra, 125, 1989,
278--287.
\bibitem [W-Y] {W-Y} K.-I. Watanabe and K. Yoshida, \emph{Hilbert-Kunz multiplicity and an inequality between multiplicity and colength}, J. Algebra, 230, 2000, 295-317.
\end{thebibliography}

\end{document}